\documentclass[12pt]{amsart}
\usepackage{amsmath}%
\usepackage{amsfonts}%
\usepackage{amssymb}%
\usepackage{graphicx}
\usepackage{tikz-cd}
\usepackage{enumerate}
\usepackage{rotating}
\usepackage{rotfloat}
\usepackage{caption}

\newtheorem{teo}{Theorem}[section]

\newtheorem{cor}[teo]{Corollary}
\newtheorem{prop}[teo]{Proposition}
\newtheorem{defi}[teo]{Definition}
\newtheorem{example}[teo]{Example}
\newtheorem{obs}[teo]{Observation}
\newtheorem{conj}[teo]{Conjecture}
\newtheorem{prob}[teo]{Problem}

\floatstyle{plain}
\floatname{diagram}{Diagram}
\newfloat{diagram}{tbp}{lop}[section]

\begin{document}

\title{Colored Tverberg theorem with new constraints on the faces \hspace{2cm}}

%\author{C. Biasi}
\author[Leandro V. Mauri]{Leandro V. Mauri}
\address{Departamento de Matem\'atica\\
	Instituto de Ci\^encias Matem\'aticas e de Computa\c c\~ao\\
	S\~ao Paulo University (USP)- C\^ampus de S\~ao Carlos \\
	13560-970, S\~ao Carlos, SP, Brazil}
\email{leandro.mauri@usp.br}

\author[D. De Mattos]{Denise de Mattos}
\address{Departamento de Matem\'atica\\
	Instituto de Ci\^encias Matem\'aticas e de Computa\c c\~ao\\
	S\~ao Paulo University (USP)- C\^ampus de S\~ao Carlos \\
	13560-970, S\~ao Carlos, SP, Brazil}
\email{deniseml@icmc.usp.br}

\author[Edivaldo L. dos Santos]{Edivaldo L. dos Santos}
\address{Departamento de Matem\'atica\\
	Universidade Federal de S\^{a}o Carlos\\
	Federal University of S\~{a}o Carlos  (UFSCAR) - C\^ampus de S\~ao Carlos \\
	13565-905, S\~ao Carlos, SP, Brazil}
\email{edivaldo@ufscar.brr}

\thanks{The first author is supported by FAPESP, process number: 2018/23928-2.}

\begin{abstract} In this paper, we prove a version of the Colored Tverberg Theorem with new constraints on the faces, in which we limit the number of faces with each one of the colors.
\end{abstract}

%To achieve proofs of Borsuk-Ulam type theorems,  several statements appears in the literature assuming X as being a cohomological n-acyclic space. In this paper, by considering a more wide class of topological spaces (not necessarily cohomological n-acyclic spaces)

\maketitle

\section{Introduction}

In 1959, Birch \cite{Bir59} formulated the following conjecture.

\vspace{0.2cm}

\begin{conj} Any $(r-1) (d+1) +1$ points in $\mathbb{R}^{d}$ can be partitioned in $N$ subsets whose convex hulls have a common point.
\end{conj}

\vspace{0.2cm}

The Birch's conjecture was proved by Helge Tverberg (see \cite{Tve66}) and since then is known as Tverberg Theorem. 

\vspace{0.2cm}

\begin{teo}[\textbf{Tverberg Theorem}] Let $d \ge 1$, $r \ge 2$ and $N= (r-1) (d+1)$ be integers. For any affine map $f: \Delta_N \rightarrow \mathbb{R}^{d}$ there are $r$ pairwise disjoint faces $\sigma_1, \ldots , \sigma_r$ of $\Delta_N$ such that $f(\sigma_1) \cap \cdots \cap f(\sigma_r) \neq \emptyset.$

\end{teo}

\vspace{0.2cm}

The following supposition is a generalization of the Tverberg Theorem to arbitrary continuous map.

\vspace{0.2cm}

\begin{conj}[\textbf{Topological Tverberg conjecture}] Let $d \ge 1$, $r \ge 2$ and $N= (r-1) (d+1)$ be integers. For any continuous map $f: \Delta_N \rightarrow \mathbb{R}^{d}$ there are $r$ pairwise disjoint faces $\sigma_1, \ldots , \sigma_r$ of $\Delta_N$ such that $f(\sigma_1) \cap \cdots \cap f(\sigma_r) \neq \emptyset.$

\end{conj}

\vspace{0.2cm}

The Topological Tverberg conjecture was considered an important unsolved problem in topological combinatorics. In 1981, the conjecture was proved, when $r$ is a prime number, by Bárány, Shlosman and Szücs \cite{Bar81}. The result for a prime number $r$ was extend for a prime power $r$, by Özaydin (unpublished) \cite{Oza87} and Volovikov \cite{Vol96}. This result is known as Topological Tverberg Theorem.

\vspace{0.2cm}

 \begin{teo}[\textbf{Topological Tverberg Theorem}] Let $d \ge 1$, $r \ge 2$ and $N= (r-1) (d+1)$ be integers. If $r$ is a prime power, then for any continuous map $f: \Delta_N \rightarrow \mathbb{R}^{d}$ there are $r$ pairwise disjoint faces $\sigma_1, \ldots , \sigma_r$ of $\Delta_N$ such that $f(\sigma_1) \cap \cdots \cap f(\sigma_r) \neq \emptyset.$

\end{teo}

\vspace{0.2cm}

The set $\{\sigma_1,\ldots ,\sigma_r \}$ of disjoint faces of $\Delta_N$ whose images by the map $f: \Delta_N \rightarrow \mathbb{R}^{d}$ have nonempty intersection is called a Tverberg partition for $f: \Delta_N \rightarrow \mathbb{R}^{d}$.

\vspace{0.2cm}

In 1992, Bárány and Larman \cite{Bar92} formulated the colored Tverberg problem, as follows.

\vspace{0.2cm}

\begin{defi}[\textbf{Coloring}] Let $N \ge 1$ be an integer and let $V(\Delta_N)$ be the set of vertices of the simplex $\Delta_N$. A \textit{coloring} of vertices $V(\Delta_N)$, by $l$ colors, is a partition $(C_1, \ldots , C_l)$ of $V(\Delta_N)$, that is, $V(\Delta_N) = C_1 \cup \cdots \cup C_l$ and $C_i \cap C_j = \emptyset$ for $1 \le i < j \le l$. The elements of the partition $(C_1, \ldots , C_l)$ are called \textit{color classes}.

\end{defi}

\vspace{0.2cm}

\begin{defi}[\textbf{Rainbow face}] Let $(C_1, \ldots , C_l)$ be the coloring of $V(\Delta_N)$ by $l$ colors. A face $\sigma$ of the  simplex $\Delta_N$ is a \textit{rainbow face} if $|\sigma \cap C_i| \le 1$, for all $1 \le i \le l$.

\end{defi}

\vspace{0.2cm}

\begin{prob}[\textbf{Bárány-Larman colored problem}] Let $d \ge 1$ and $r \ge 2$ be integers. Determine the smallest number $n= n(d,r)$ such that for every map $f: \Delta_{n-1} \rightarrow \mathbb{R}^{d}$ and every coloring $(C_1, \ldots , C_{d+1})$ of the vertex set $V(\Delta_{n-1})$ of the simplex $\Delta_{n-1}$ by $d+1$ colors, with each color of size at least $r$, there exist $r$ pairwise disjoint faces $\sigma_1, \ldots , \sigma_r$ of $\Delta_{n-1}$ satisfying \footnote{$[n]=\{1, \cdots,n\}.$}: \begin{eqnarray} {\rm (i)} \hspace{0.2cm} | C_i \cap \sigma_j| \le 1 \hspace{0.1cm} \mbox{, for every } i \in [d+1], \hspace{0.1cm} j \in [r], \hspace{0.1cm} \mbox{that is, } \sigma_1, \ldots, \sigma_r \hspace{0.1cm} \nonumber \\ \mbox{are rainbow faces; and} \nonumber \hspace{7.2cm} \\ {\rm (ii)} \hspace{0.2cm} f(\sigma_1) \cap \cdots \cap f(\sigma_r) \neq \emptyset. \hspace{6.7cm}\nonumber\end{eqnarray}

\end{prob}

\vspace{0.2cm}

A modified colored Tverberg problem was presented by \v{Z}ivaljevi\'{c} and Vre\'{c}ica in the paper \cite{Ziv92}.

\vspace{0.2cm}

\begin{prob}[\textbf{The \v{Z}ivaljevi\'{c} and Vre\'{c}ica colored Tverberg problem}] Let $d \ge 1$ and $r \ge 2$ be integers. Determine the smallest number $t=t(d,r)$ such that for every affine (or continuous) map $f: \Delta \rightarrow \mathbb{R}^{d}$, and every coloring $(C_1, \ldots , C_{d+1})$ of the the vertex set $V(\Delta)$ by $d+1$ colors, with each color of size at least $t$, there exist $r$ pairwise disjoint faces $\sigma_1, \ldots , \sigma_r$ of $\Delta$ satisfying: \begin{eqnarray} {\rm (i)} \hspace{0.2cm} | C_i \cap \sigma_j| \le 1 \hspace{0.1cm} \mbox{, for every } i \in [d+1], \hspace{0.1cm} j \in [r], \hspace{0.1cm} \mbox{that is, } \sigma_1, \ldots, \sigma_r \hspace{1.3cm} \nonumber\\ \mbox{are rainbow faces; and} \nonumber \hspace{8.4cm} \\ {\rm(ii)} \hspace{0.2cm} f(\sigma_1) \cap \cdots \cap f(\sigma_r) \neq \emptyset. \hspace{8cm}\nonumber\end{eqnarray}

\end{prob}

\vspace{0.2cm}

For $r \ge 2$ a prime power, \v{Z}ivaljevi\'{c} and Vre\'{c}ica proved that $t(d,r) \le 2r-1$. This result was known as colored Tverber theorem of \v{Z}ivaljevi\'{c} and Vre\'{c}ica.

\vspace{0.2cm}

\begin{teo}[\textbf{Colored Tverberg theorem of \v{Z}ivaljevi\'{c} and Vre\'{c}ica} \cite{Ziv92}] Let $d \ge 1$ be an integer and let $r \ge 2$ be a prime power. For every continuous map $f: \Delta \rightarrow \mathbb{R}^{d}$ and every coloring $(C_1, \ldots , C_{d+1})$ of the the vertex set $V(\Delta)$ by $d+1$ colors, with each color of size at least $2r-1$, there exist $r$ pairwise disjoint faces $\sigma_1, \ldots , \sigma_r$ of $\Delta$ satisfying: \begin{eqnarray} {\rm(i)} \hspace{0.2cm} | C_i \cap \sigma_j| \le 1 \hspace{0.1cm} \mbox{, for every } i \in [d+1], \hspace{0.1cm} j \in [r], \hspace{0.1cm} \mbox{that is, } \sigma_1, \ldots, \sigma_r \hspace{1.7cm} \nonumber \\ \mbox{are rainbow faces; and} \hspace{8.7cm} \nonumber \\ {\rm (ii)} \hspace{0.2cm} f(\sigma_1) \cap \cdots \cap f(\sigma_r) \neq \emptyset. \hspace{8.4cm}\nonumber\end{eqnarray}

\vspace{0.2cm}

Starting from Colored Tverberg theorem of \v{Z}ivaljevi\'{c} and Vre\'{c}ica, we introduce a problem that increases the restrictions on the faces.

\vspace{0.2cm}

\begin{prob} \label{Prob1.9}

Let $d \ge 1$, $m\ge 1$ be integers and let $r \ge 2$ be a prime power. Determine integers $ l_1, \ldots , l_m $ $( 1 \le l_i \le r$, for $i=1, \cdots , m)$ such that for every continuous map $f: \Delta \rightarrow \mathbb{R}^{d}$ and every coloring $(C_1, \ldots , C_{m})$ of the the vertex set $V(\Delta)$ by $m$ colors, with each color of size at least $2r-1$, there exist $r$ pairwise disjoint faces $\sigma_1, \ldots , \sigma_r$ of $\Delta$ satisfying: \begin{eqnarray} {\rm (i)} \hspace{0.2cm} | C_i \cap \sigma_j| \le 1 \hspace{0.1cm} \mbox{, for every } i \in [m], \hspace{0.1cm} j \in [r], \hspace{0.1cm} \mbox{that is, } \sigma_1, \ldots, \sigma_r \hspace{1.7cm} \nonumber \\ \mbox{are rainbow faces} ;\nonumber \hspace{8.8cm} \\ {\rm (ii)} \hspace{0.2cm} f(\sigma_1) \cap \cdots \cap f(\sigma_r) \neq \emptyset \mbox{; and} \hspace{7cm} \nonumber \\ {\rm (iii)} \hspace{0.2cm} | (\sigma_1 \cup \cdots \cup \sigma_r) \cap C_i| \le l_i \hspace{0.1cm} \mbox{, for every } i \in [m]. \hspace{4cm} \nonumber \end{eqnarray}

\end{prob}

\vspace{0.2cm}

In this paper, we prove that the Problem $1.9$ is true, when \begin{eqnarray}  \displaystyle \sum_{i=1}^{m} l_i > (d+1)(r-1) \hspace{0.1cm} \mbox{(see Theorem 5.1)}. \nonumber \end{eqnarray}

\vspace{0.2cm}

\section{The index of Volovikov}

\vspace{0.2cm}

In this section, we present the index of Volovikov, which is the algebraic topology tool that we use to prove our main results in Section $5$. For more details on the index of Volovikov, see \cite{Vol00}.

\vspace{0.2cm}

Initially, we define the \textit{Borel construcion}.

\vspace{0.2cm}

\begin{defi}[\textbf{Borel construction}] Let $G$ be a compact Lie group and let $X$ be a Hausdorff, paracompact  $G$-space, on which the $G$ acts freely. Then $X \rightarrow X / G$ is a principal $G$-bundle (see \cite{Bre72}) and one can take  \begin{eqnarray} h: X / G \rightarrow BG\nonumber\end{eqnarray} a classifying map for the $G$-bundle $X \rightarrow X / G$.

\vspace{0.2cm}

Let us consider the product $EG \times X$, with diagonal $G$-action $g(e,x) =(ge, gx)$, for every $g \in G$ and $(e,x) \in EG \times X$. Let $EG \times_G X = (EG \times X) / G$ be its orbit space. The first projection $EG \times X \rightarrow EG$ induces a map: \begin{eqnarray} p_X: EG \times_G X \rightarrow (EG)/G = BG, \nonumber \\ (e,x) G \mapsto e G  \hspace{1.9cm} \nonumber \end{eqnarray} which is a fibration with fiber $X$ and base space $BG$ being the classifying space of $X$.

\vspace{0.2cm}

This is the \textit{Borel construction}. It associates to each Hausdorff, paracompact $G$-space $X$, a $G$-space $EG \times_G X$, which is denoted by $X_G$, over $BG$. Also, it associates to each $G$-map $f: X \rightarrow Y$, a fiber preserving map $\overline{ \mbox{Id}_{EG} \times f}:EG \times_G X \rightarrow EG \times_G Y$.

\end{defi}

\vspace{0.2cm}

Let us recall the Leray-Serre theorem for fibrations.

\vspace{0.2cm}

\begin{teo}[\textbf{The cohomology Leray-Serre Spectral sequence (Theorem $5.2$ \cite{McC01})}]
Let $R$ be a commutative ring with unit. Given a fibration $F \hookrightarrow E \stackrel{p}{\rightarrow} B$, where $B$ is a path-wise connected space, there is a first quadrant spectral sequence of algebras $\{ E_r^{\ast, \ast}, d_r \}$, with: \begin{eqnarray} E_{2}^{p,q} \cong H^{p} ( B ; \mathcal{H}^{q} (F;R)), \nonumber \end{eqnarray} the cohomology of $B$, with local coefficients in the cohomology of $F$, the fiber of $p$, and converging to $H^{\ast}(E;R)$ as an algebra. Furthermore, this spectral sequence is natural with respect to fiber-preserving maps of fibrations.

\end{teo}

\vspace{0.2cm}

Let us now recall one of the numerical index defined by Volovikov in \cite{Vol00}. This is a function on $G$-spaces, whose value is either a positive integer or $\infty$. For our purposes, it is sufficient to consider that $G$ is a $p$-torus ($p$ a prime number), that is, $G= (\mathbb{Z}_p)^{n}$, for $n \ge 1$ and $\mathbb{K} = \mathbb{Z}_p$ is a coefficient field in the \v{C}ech cohomology .

\vspace{0.2cm}

\begin{defi}[\textbf{The index of Volovikov}] Let $G= (\mathbb{Z}_p)^{n}$ ($p$ a prime number) be a compact Lie group, for $n \ge 1$ and let $X$ be a Hausdorff paracompact $G$-space, on which the $G$ acts freely. The definition of \textit{the Volovikov index of $X$}, denoted by $i(X)$, uses the spectral sequence of the bundle $p_X: X_G \rightarrow BG$, with fibre $X$ (the Borel construction), given in Theorem $2.2$. This spectral sequence this of a point converges to the equivariant cohomology $H^{\ast} (X_G; \mathbb{Z}_p)$. Let $\Lambda^{\ast}$ be the equivariant cohomology algebra $H^{\ast} (\mbox{pt}_G; \mathbb{Z}_p)= H^{\ast} (BG; \mathbb{Z}_p)$. Suppose that $X$ is path connected. Then $E_{2}^{\ast,0} = \Lambda^{\ast}$. Assume that $E_{2}^{\ast,0}= \cdots = E_{s}^{\ast,0} \neq E_{s+1}^{\ast,0}$. Then, by definition, $i(X)=s$. If $E_{2}^{\ast,0}= \cdots = E_{\infty}^{\ast,0}$ then, by definition, $i(X)= \infty$.

\end{defi}

\vspace{0.2cm}

We state some properties of the index of Volovikov (see \cite{Vol00}).

\vspace{0.2cm}

\begin{prop} Let $X$, $Y$ and $Z$ be a Hausdorff, paracompact  $G$-spaces, on which the $G$ acts freely. Then:

\vspace{0.2cm}

{\rm (i)} \vspace{0.1cm} If there is a $G$-equivariant map $X \rightarrow Y$, then $i(X) \le i(Y)$.

{\rm (ii)} \vspace{0.1cm} If $\tilde{H}^{j}(X; \mathbb{Z}_p)=0$, for all $j<n$, then $i(X) \ge n+1$.

{\rm (iii)} \vspace{0.1cm} If $H^{j}(Z; \mathbb{Z}_p)=0$, for all $j \ge n$ and if $i(Z)< \infty$, then $i(Z) \le n$.

{\rm (iv)} \vspace{0.1cm} If $X$ is a compact or finite dimensional space such that $H^{\ast}(X; \mathbb{Z}_p) = H^{\ast}(S^{n}; \mathbb{Z}_p) $, and if $G$ acts on $X$ without fixed points, then $i(X) =n+1$.

\end{prop}

\vspace{0.2cm}

\section{Chessboard complex and connectedness}

\vspace{0.2cm}

In this section, we introduce the chessboard complex, which is widely used in Colored Tverberg Theorem proofs. Also, we introduce the concept of connectedness and we relate it to the index of Volovikov. For more details on these topics, see \cite{Bla17} and \cite{mat08}.

\vspace{0.2cm}

\begin{defi}[\textbf{Chessboard complex}] The \textit{$m \times n$ chessboard complex} $\Delta_{m,n}$ is the simplicial complex whose vertex set is $[m] \times [n]$. The simplexes of $\Delta_{m,n}$ are the subsets $\{(i_0,j_0), \ldots, (i_k,j_k) \} \subset [m] \times [n]$, where $i_s \neq i_{s'}$ ($1 \le s < s'\le k$), and $j_t \neq j_{t'}$ ($1 \le t < t'\le k$).

\end{defi}

\vspace{0.2cm}

\begin{defi} Let $n \ge -1$ be an integer. A topological space $X$ is \textit{$n$-connected} if any continuous map $f: S^{k} \rightarrow X$, where $-1 \le k \le n$, can be continuously extended to a continuous map $g: B^{k+1} \rightarrow X$, that is, $g|_{\partial B^{k+1} = S^k} = f$ (here $B^{k+1}$ denotes a $(k+1)$-dimensional closed ball whose boundary is the sphere $S^{k}$). A topological space is $(-1)$-connected if it is non-empty. If the space $X$ is $n$-connected, but it is not $(n+1)$-connected, we write $\mbox{conn} (X)=n$.

\end{defi}

\vspace{0.2cm}

\begin{teo}[pg. 332, \cite{Bla17}] Let $X$ and $Y$ be topological spaces. Then: \begin{eqnarray} \mbox{conn} (X \ast Y) \ge  \mbox{conn} (X) + \mbox{conn} (Y) +2. \nonumber \end{eqnarray} \end{teo}

\vspace{0.2cm}

\begin{teo}[4.4.1 Theorem \cite{mat08}] Let $X$ be a nonempty topological space and let $k \ge 1$ be an integer. Then, $X$ is $k$-connected, if and only if, it is simply connected (i.e., the fundamental group $\pi_1(X)$ is trivial) and $\tilde{H}_i(X)=0$, for all $i=0,1, \ldots,k$.

\end{teo}

\vspace{0.2cm}

\begin{teo} Let $X$ be a topological space. Then, $i (X ) \ge  \mbox{conn} (X) +2.$
\end{teo}

\vspace{0.2cm}

\textit{Proof.} It is a consequence of Theorem $3.4$ and Proposition $2.4 $ ${\rm (ii)}$.

\vspace{0.1cm}

\hspace{12.3cm} $\square$

\vspace{0.2cm}

\begin{teo} [4.4.2 Proposition \cite{mat08}] A simplicial complex $K$ is $s$-connected, if and only if, the $(s+1)$-skeleton $K^{\le s+1}$ is $s$-connected.
\end{teo}

\vspace{0.2cm}

\begin{teo} Let $m, n \ge 1$ be integers. Then, \begin{eqnarray} \mbox{conn} (\Delta_{m,n}) = \mbox{min} \left\{m,n, \left\lfloor{\frac{m+n+1}{3}}\right\rfloor \right\}-2. \nonumber \end{eqnarray}

\end{teo}

\vspace{0.2cm}

\begin{cor} Let $r \ge 2$ be an integer. Then,  $\mbox{conn} (\Delta_{2r-1,r}) = r-2.$

\end{cor}

\vspace{0.2cm}

\section{Key Theorem}

\vspace{0.2cm}

In this section, we are going to adapt Corollary $4.8$ of the paper \cite{Bla17}, which is fundamental for the proof of  the Colored Tverberg Theorem of \v{Z}ivaljevi\'{c} and Vre\'{c}ica. We make the necessary adjustment which include the new constraint \textbf{(iii)} of the Problem $1.9$ on the faces.

\vspace{0.2cm}

Let us observe that the \textit{deleted join} $(\Delta_N)_{\Delta(2)}^{\ast r} = \{ \sigma_1 \uplus \cdots \uplus \sigma_r \hspace{0.1cm} \mbox{;} \hspace{0.1cm} \sigma_i \cap \sigma_j = \emptyset, \mbox{for } i,j \in [r] \hspace{0.1cm} \mbox{and } i \neq j \}$ can be seen as the join of the collections of $r$ pairwise disjoint faces of $\Delta_N$.

\vspace{0.2cm}

If we have a coloring $(C_1, \ldots , C_m)$ of the vertex set $V(\Delta)$ by $m$ colors, the \textit{$r$- fold $2$-wise deleted join} of rainbow subcomplex $R_{(C_1, \cdots , C_m)}$ can be identified with: \begin{eqnarray} \left( R_{(C_1, \ldots , C_m)}\right)_{\Delta(2)}^{\ast r} \cong \left( [|C_1|] \ast \cdots \ast [|C_m|] \right)_{\Delta(2)}^{\ast r} \cong \Delta_{|C_1|,r} \ast \cdots \ast \Delta_{|C_m|,r}. \nonumber \end{eqnarray}

\vspace{0.2cm}

The action of the symmetric group $\mathfrak S_{r} = \mbox{Sym} (r)$ on the chessboard complex $\Delta_{|C_i|,r}$, for $i=1, \cdots, m$, is given by the permutation of columns of the chessboard, that is \begin{eqnarray}  \pi \cdot \{ (i_0, j_0) , \ldots , (i_k, j_k) \} = \{ (i_0, \pi(j_0)) , \ldots , (i_k, \pi(j_k)) \},  \nonumber \\ \hspace{0.1cm} \mbox{for } \pi \in \mathfrak S_{r} \hspace{0.1cm} \mbox{and} \hspace{0.1cm} \{ (i_0, j_0) , \ldots , (i_k, j_k) \} \in \Delta_{|C_i|,r}. \hspace{1cm}  \nonumber \end{eqnarray}

\vspace{0.2cm}

Then, $\Delta_{|C_1|,r} \ast \cdots \ast \Delta_{|C_m|,r}$ can be seen as the join of the collections of $r$ pairwise disjoint rainbow faces, that is, it satisfies the condition ${\rm (i)}$.

\vspace{0.2cm}

To obtain a join of $r$ disjoint faces satisfying the conditions ${\rm(i)}$ and ${\rm (iii)}$ of the Problem $1.9$, that is, the faces are rainbow and with at most $l_i$ vertices of color $i$, for $i=1, \ldots , m$, we just need to take the $(l_i-1)$-skeleton of the chessboard $\Delta_{|C_i|,r}$. Therefore, the join \begin{eqnarray}\Delta_{|C_1|,r}^{ \le l_1-1} \ast \cdots \ast \Delta_{|C_m|,r}^{\le l_m-1} \nonumber \end{eqnarray} can be seen as the join of the collections of $r$ pairwise disjoint faces satisfying the conditions {\rm (i)} and {\rm (iii)} of the Problem $1.9$.

\vspace{0.2cm}

\begin{teo} Let $d \ge 1$, $m\ge 1$, $l_1, \ldots, l_m$ ($1 \le l_i \le r$, for $i=1, \ldots , m$) be integers and let $r =p^{n} \ge 2$ be a prime power. Let $(C_1, \ldots , C_m)$ be a coloring of $\Delta$ by $m$ colors. If there is no a $\mathfrak S_r$-equivariant map \footnote{$W_r= \{ (t_1, \ldots, t_r) \in \mathbb{R}^{r} \hspace{0.1cm}; \hspace{0.1cm} \sum_{i=1}^{r} t_i =0\}.$} \begin{eqnarray} \Delta_{|C_1|,r}^{\le l_1-1} \ast \cdots \ast \Delta_{|C_m|,r}^{ \le l_m-1}  \longrightarrow S \left(W_{r}^{\oplus (d+1)}\right)\nonumber \end{eqnarray} then for every continuous map $f: \Delta \rightarrow \mathbb{R}^{d}$, there exist $r$ pairwise disjoint faces $\sigma_1, \ldots, \sigma_r$ of $\Delta$ satisfying \begin{eqnarray} {\rm (i)} \hspace{0.2cm} | C_i \cap \sigma_j| \le 1 \mbox{,} \hspace{0.1cm} \mbox{ for every } i \in [m], \hspace{0.1cm} j \in [r], \hspace{0.1cm} \mbox{that is, } \sigma_1, \ldots, \sigma_r \hspace{2cm} \nonumber \\ \mbox{are rainbow faces} ;\nonumber \hspace{9.1cm} \\ {\rm(ii)} \hspace{0.2cm} f(\sigma_1) \cap \cdots \cap f(\sigma_r) \neq \emptyset \mbox{; and} \hspace{7.4cm} \nonumber \\ {\rm (iii)} \hspace{0.2cm} | (\sigma_1 \cup \cdots \cup \sigma_r) \cap C_i| \le l_i \mbox{,} \hspace{0.1cm} \mbox{ for every } i \in [m]. \hspace{4.4cm} \nonumber \end{eqnarray}

\end{teo}

\vspace{0.2cm}

\textit{Proof.} Let us consider $f: \Delta \rightarrow \mathbb{R}^{d}$ a continuous map such that for any $r$ pairwise disjoint faces $\sigma_1, \ldots, \sigma_r$ of $\Delta$ which satisfy the conditions ${\rm (i)}$ and ${\rm (iii)}$, the constraint $\textbf{(ii)}$ is not satisfied, that is $f(\sigma_1) \cap \cdots \cap f(\sigma_r) = \emptyset$. That is, $f$ is a counterexample for the statement of the theorem.

\vspace{0.2cm}

Let us define the \textit{join map} \begin{eqnarray} J_f : \left( \Delta_N \right)_{\Delta(2)}^{\ast r} \longrightarrow \left(\mathbb{R}^{d+1}\right)^{\oplus r}\nonumber \hspace{4cm}\\ \lambda_1 x_1 + \cdots +\lambda_r x_r \longmapsto (\lambda_1, \lambda_1 f(x_1)) \oplus \cdots \oplus ( \lambda_r, \lambda_r f(x_r)).\nonumber\end{eqnarray}

\vspace{0.2cm}

Note that $\Delta_{|C_1|,r}^{\le l_1 -1} \ast \cdots \ast \Delta_{|C_m|,r}^{\le l_m -1} $ is a $\mathfrak S_r$-equivariant subespace of $(\Delta_N)_{\Delta(2)}^{\ast r}$, since the permutation of the $r$ faces does not change the number of vertices of color $i$, for each $i \in [m]$. 

\vspace{0.2cm}

Let us define the \textit{join map} \begin{eqnarray} J'_f \stackrel{def}{=} J_f|_{\Delta_{|C_1|,r}^{\le l_1 -1} \ast \cdots \ast \Delta_{|C_m|,r}^{\le l_m -1}} : \Delta_{|C_1|,r}^{\le l_1 -1} \ast \cdots \ast \Delta_{|C_m|,r}^{\le l_m -1} \longrightarrow \left(\mathbb{R}^{d+1}\right)^{\oplus r}\nonumber \hspace{4cm}\\  \lambda_1 x_1 + \cdots +\lambda_r x_r \longmapsto (\lambda_1, \lambda_1 f(x_1)) \oplus \cdots \oplus ( \lambda_r, \lambda_r f(x_r)). \hspace{5cm} \nonumber\end{eqnarray}

\vspace{0.2cm}

Let us consider the $\mathfrak S_r$- invariant subespace \begin{eqnarray} D_J= \{ (z_1, \ldots , z_r) \in (\mathbb{R}^{d+1})^{\oplus r} \hspace{0.1cm} \mbox{;} \hspace{0.1cm} z_1 = \cdots = z_r \}. \nonumber \end{eqnarray}

\vspace{0.2cm}

We have that $f$ is a counterexample which does not satisfy the condition ${\rm (ii)}$, therefore $\mbox{im} (J'_f) \cap D_J = \emptyset$.

\vspace{0.2cm}

Indeed, let us assume that \begin{eqnarray} (\lambda_1, \lambda_1 f(x_1)) \oplus \cdots \oplus ( \lambda_r, \lambda_r f(x_r)) \in \mbox{im} (J'_f) \cap D_J \neq \emptyset \mbox{,}\nonumber \end{eqnarray} for some $\lambda_1 x_1 + \cdots +\lambda_r x_r \in \Delta_{|C_1|,r}^{\le l_1 -1} \ast \cdots \ast \Delta_{|C_m|,r}^{\le l_m -1}$. Then, $\lambda_1= \cdots = \lambda_r = \frac{1}{r}$ and, consequently, $f(x_1)= \cdots = f(x_r)$, where $ x_1 \in \mbox{relint } \sigma_1, \ldots , x_r \in \mbox{relint } \sigma_r$ \footnote{The relative interior of $\sigma$, denoted by $\mbox{relint } \sigma$, in which the vertex set of $\sigma$ is $V(\sigma)= \{v_0, \ldots, v_{d+1}\}$, is given by $\mbox{relint } \sigma = \{ \lambda_0 v_0 + \cdots + \lambda_{d+1} v_{d+1} \hspace{0.1cm} ; \hspace{0.1cm}  0 < \lambda_0, \ldots , \lambda_{d+1} \le 1 \hspace{0.1cm} \mbox{and} \hspace{0.1cm} \sum_{i=0}^{d+1} \lambda_i = 1 \}$. } and these $r$ pairwise disjoint faces satisfy the conditions {\rm (i)} and {\rm (iii)}, as already explained previously. This contradicts the fact that $f$ is a counterexample for the theorem.

\vspace{0.2cm}

Thus, $J'_f$ induces an $\mathfrak S_r$-equivariant map \begin{eqnarray}  \Delta_{|C_1|,r}^{\le l_1 -1} \ast \cdots \ast \Delta_{|C_m|,r}^{\le l_m -1} \longrightarrow \left(\mathbb{R}^{d+1}\right)^{\oplus r} \setminus D_J \mbox{,}\nonumber \end{eqnarray}
which we also denote by $J'_f$, with a slight of abuse of notation.

\vspace{0.2cm}

Furthermore, let \begin{eqnarray} R_J: \left(\mathbb{R}^{d+1}\right)^{\oplus r} \setminus D_J \rightarrow D_J^{\perp} \setminus \{0\} \rightarrow S\left(D_J^{\perp}\right)\nonumber \end{eqnarray} be the composition of the appropiate projection and deformation retraction. The map $R_J$ is a $\mathfrak S_r$-equivariant.

\vspace{0.2cm}

We have an isomorphism of $\mathfrak S_r$-representations $D_J^{\perp} \cong W_r^{\oplus (d+1)}$, where $W_r =\{( t_1, \cdots, t_r) \in \mathbb{R}^{r} \hspace{0.1cm} \mbox{;} \hspace{0.1cm} \sum_{i=1}^{r} t_i =0 \}$ and it is equipped with the action of the symmetric group $\mathfrak S_r$ given by permutation of coordinates.

\vspace{0.2cm}

Identifying $D_J^{\perp}$ with $ W_r^{\oplus (d+1)}$ we have that $R_J$ is a $\mathfrak S_r$-equivariant map \begin{eqnarray} R_J: \left(\mathbb{R}^{d+1}\right)^{\oplus r} \setminus D_J \rightarrow S\left(W_{r}^{\oplus(d+1)}\right).\nonumber \end{eqnarray}

\vspace{0.2cm}

Then, we have a $\mathfrak S_r$-equivariant map, given by the composition \begin{eqnarray} R_J \circ J'_f : \Delta_{|C_1|,r}^{\le l_1 -1} \ast \cdots \ast \Delta_{|C_m|,r}^{\le l_m -1}  \longrightarrow S \left(W_{r}^{\oplus (d+1)}\right). \nonumber \end{eqnarray}

\vspace{0.2cm}

The existence of the map $R_J \circ J'_f$ is a contradiction to the hypothesis, so the theorem is proved.

\hspace{12.3cm} $\square$

\vspace{0.2cm}

\section{Colored Tverberg Theorem with new constraints on the faces}

\vspace{0.2cm}

In this section, we enunciate and demonstrate the new results.

\vspace{0.2cm}

\begin{teo}[\textbf{Colored Tverberg Theorem with new constraint on the faces}]

Let $d \ge 1$ be an integer, let $r = p^{n}\ge 2$ be a prime power and let $ l_1, \ldots , l_m $ $( 1 \le l_i \le r$, for $i=1, \ldots , m)$ be  integers satisfying $ \displaystyle \sum_{i=1}^{m} l_i > (d+1)(r-1)$ . For every continuous map $f: \Delta \rightarrow \mathbb{R}^{d}$ and for every coloring $(C_1, \ldots , C_{m})$ of the vertex set of the simplex $\Delta$ by $m$ colors, with each color of size at least $2r-1$, there exist $r$ pairwise disjoint faces $\sigma_1, \ldots , \sigma_r$ of $\Delta$  satisfying: \begin{eqnarray} {\rm (i)} \hspace{0.2cm} | C_i \cap \sigma_j| \le 1 \mbox{,} \hspace{0.1cm} \mbox{ for every } i \in [m],  j \in [r], \hspace{0.1cm} \mbox{that is, } \sigma_1, \ldots, \sigma_r \hspace{2.1cm} \nonumber \\ \mbox{are rainbow faces} ;\nonumber \hspace{9cm} \\ {\rm (ii)} \hspace{0.2cm} f(\sigma_1) \cap \cdots \cap f(\sigma_r) \neq \emptyset \nonumber \mbox{; and} \hspace{7.3cm} \\ {\rm (iii)} \hspace{0.2cm} | (\sigma_1 \cup \cdots \cup \sigma_r) \cap C_i| \le l_i \mbox{,} \hspace{0.1cm} \mbox{ for every } i \in [m]. \hspace{0.1cm}  \hspace{4.2cm}\nonumber \end{eqnarray}

\end{teo}

\vspace{0.2cm}

\textit{Proof.} Without loss of generality, we can assume that $|C_1| = \cdots = |C_m| = 2r-1$. Then, by Theorem $4.1$, it is sufficient to show that there is no a $\mathfrak S_r$-equivariant map \begin{eqnarray}\Delta_{2r-1,r}^{\le l_1 -1} \ast \cdots \ast \Delta_{2r-1,r}^{\le l_m -1} \longrightarrow S \left(W_{r}^{\oplus (d+1)}\right). \nonumber \end{eqnarray}

\vspace{0.2cm}

Consider the regular embedding \begin{eqnarray} \varphi: (\mathbb{Z}_p)^{n} \rightarrow \mbox{Sym } \left((\mathbb{Z}_p)^{n}\right) \cong \mathfrak S_r  \nonumber \\ g \longmapsto L_g : (\mathbb{Z}_p)^{n} \rightarrow (\mathbb{Z}_p)^{n}, \hspace{0.1cm} L_g (x) =g+x.  \nonumber\end{eqnarray}

\vspace{0.2cm}

Therefore, we have a subgroup $(\mathbb{Z}_p)^{n} \cong \mbox{Im} (\varphi) \le \mathfrak S_r$. Then, to prove the non-existence of a $\mathfrak S_r$-equivariant map it suffices to prove the non-existence of a  $(\mathbb{Z}_p)^{n}$-equivariant map \begin{eqnarray}\Delta_{2r-1,r}^{\le l_1 -1} \ast \cdots \ast \Delta_{2r-1,r}^{\le l_m -1} \longrightarrow S \left(W_{r}^{\oplus (d+1)}\right). \nonumber \end{eqnarray}

\vspace{0.2cm}

Let us suppose that there is a such $(\mathbb{Z}_p)^{n}$-equivariant map. Then, by Proposition $2.4 (i)$, we have that \begin{eqnarray}i\left(\Delta_{2r-1,r}^{\le l_1 -1} \ast \cdots \ast \Delta_{2r-1,r}^{\le l_m -1}\right) \le i\left( S \left(W_{r}^{\oplus (d+1)}\right) \right). \nonumber \end{eqnarray}

\vspace{0.2cm}

By Proposition $2.4(iv)$, we have that  \begin{eqnarray} i\left( S \left(W_{r}^{\oplus (d+1)}\right) \right) = i\left( S^{(r-1)(d+1)-1}\right) = (r-1) (d+1). \nonumber \end{eqnarray}

\vspace{0.2cm}

On the other hand, note that $\Delta_{2r-1,r}^{\le l_i -1}$ is $(l_i-2)$-connected, for any $i=1, \ldots , m$.

\vspace{0.2cm}

Indeed, $\Delta_{2r-1,r}$ is $(r-2)$-connected (by Corollary $3.8$) so, it is $(l_i-2)$-connected (since $l_i \le r$). Hence, by Theorem $3.6$, we have that $\Delta_{2r-1,r}^{\le l_i -1}$ is $(l_i-2)$-connected.

\vspace{0.2cm}

By Theorem $3.3$, we have that \begin{eqnarray} \mbox{conn} \left( \Delta_{2r-1,r}^{\le l_1 -1} \ast \cdots \ast \Delta_{2r-1,r}^{\le l_m -1} \right) \ge \left( \displaystyle \sum_{i=1}^{m} (l_i -2) \right) + 2 (m-1) = \left( \displaystyle \sum_{i=1}^{m} l_i \right)-2.  \nonumber \end{eqnarray}

\vspace{0.2cm}

It follows from, by Theorem $3.5$ and the inequality below \begin{eqnarray} \displaystyle \sum_{i=1}^{m} l_i > (d+1)(r-1) \nonumber \end{eqnarray} , that \begin{eqnarray} i \left(\Delta_{2r-1,r}^{\le l_1 -1} \ast \cdots \ast \Delta_{2r-1,r}^{\le l_m -1}\right) \ge \left[  \left( \displaystyle \sum_{i=1}^{m} l_i  \right)-2 \right]+2 = \displaystyle \sum_{i=1}^{m} l_i > (r-1) (d+1). \nonumber \end{eqnarray}

\vspace{0.2cm}

Then, $i\left(\Delta_{2r-1,r}^{\le l_1 -1} \ast \cdots \ast \Delta_{2r-1,r}^{\le l_m -1}\right) > (r-1) (d+1) =i \left(S\left(W_{r}^{\oplus(d+1)}\right)\right)$, but it contradicts \begin{eqnarray}i\left(\Delta_{2r-1,r}^{\le l_1 -1} \ast \cdots \ast \Delta_{2r-1,r}^{\le l_m -1}\right) \le i \left(S\left(W_{r}^{\oplus(d+1)}\right)\right). \nonumber \end{eqnarray} 

\vspace{0.2cm}

Therefore, there is no a $(\mathbb{Z}_p)^n$-equivariant map  \begin{eqnarray} \Delta_{2r-1,r}^{\le l_1 -1} \ast \cdots \ast \Delta_{2r-1,r}^{\le l_m -1} \longrightarrow S \left(W_{r}^{\oplus (d+1)}\right).  \nonumber \end{eqnarray}
 
 \hspace{12.3cm}$\square$

\vspace{0.2cm}

\begin{obs} Theorem $5.1$ provides an answer to the Problem $1.9$ proposed previously.
\end{obs} 

\vspace{0.2cm}

\begin{cor}[\textbf{Colored Tverberg theorem of \v{Z}ivaljevi\'{c} and Vre\'{c}ica with new constraints on the faces}]

Let $d \ge 1$ be an integer, let $r =p^{n}\ge 2$ be a prime power and let  $ l_1, \cdots , l_{d+1} $ $( 1 \le l_i \le r$, for $i=1, \cdots , d+1)$ be  integers satisfying $ \displaystyle \sum_{i=1}^{d+1} l_i > (d+1)(r-1)$ . For every continuous map $f: \Delta \rightarrow \mathbb{R}^{d}$ and for every coloring $(C_1, \ldots , C_{d+1})$ of the vertex set of the simplex $\Delta$ by $d+1$ colors, with each color of size at least $2r-1$, there exist $r$ pairwise disjoint faces $\sigma_1, \ldots , \sigma_r$ of $\Delta$  satisfying \begin{eqnarray} {\rm (i)} \hspace{0.2cm} | C_i \cap \sigma_j| \le 1 \mbox{,} \hspace{0.1cm} \mbox{ for every } i \in [d+1], \hspace{0.1cm} j \in [r], \hspace{0.1cm} \mbox{that is, } \sigma_1, \ldots, \sigma_r \hspace{2.1cm} \nonumber \\ \mbox{are rainbow faces;} \hspace{10cm} \nonumber \\ {\rm (ii)} \hspace{0.2cm} f(\sigma_1) \cap \cdots \cap f(\sigma_r) \neq \emptyset \hspace{0.1cm} \mbox{; and}  \hspace{7.9cm}\nonumber \\ {\rm (iii)} \hspace{0.2cm} | (\sigma_1 \cup \cdots \cup \sigma_r) \cap C_i| \le l_i \mbox{,} \hspace{0.1cm} \mbox{ for every } i \in [d+1]. \hspace{4.4cm}\nonumber\end{eqnarray}

\end{cor}

\vspace{0.2cm}

\textit{Proof.} Just apply Theorem $5.1$ to the case that $m=d+1$.

\hspace{12.3cm} $\square$

\vspace{0.2cm}

\begin{obs} The Corollary $5.3$ gives a new version of the Colored Tverberg theorem of \v{Z}ivaljevi\'{c} and Vre\'{c}ica  with more constraints on the $r$ faces.

\vspace{0.2cm}

The following example shows how Corollary $5.3$ is a stronger version of the Colored Tverberg theorem of \v{Z}ivaljevi\'{c} and Vre\'{c}ica.

\vspace{0.2cm}

\begin{example}

Let $d=2$, $r=3$ and let $\mathcal{C} = (C_1, C_2, C_3)$ be a coloring of the vertex set $V(\Delta)=[15]$, with $C_1=\{1,2,3,4,5\}$, $C_2=\{6,7,8,9,10\}$ and $C_3 = \{11,12,13,14,15\}$, that is, $|C_1|=|C_2|=|C_3| =2r-1=5$. Let $f: \Delta \rightarrow \mathbb{R}^{2}$ be a continuous map. 

\vspace{0.2cm}

By Colored Tverberg theorem of \v{Z}ivaljevi\'{c} and Vre\'{c}ica (Theorem $1.8$), we can have a Tverberg partition $\sigma_1=\{1,6,12\}$, $\sigma_2=\{5,9,14\}$ and $\sigma_3= \{10,11\}$ such that $f(\sigma_1) \cap f(\sigma_2) \cap f(\sigma_3) \neq \emptyset$.

\vspace{0.2cm}

In Corollary $5.3$, with $l_1=l_2=2$ and $l_3=3$, the Tverberg partition $\sigma_1=\{1,6,12\}$, $\sigma_2=\{5,9,14\}$ and $\sigma_3= \{10,11\}$ does not satisfy  the condition {\rm (iii)} in the Problem $1.9$, since $|(\sigma_1 \cup \sigma_2 \cup \sigma_3) \cap C_2| = 3 > l_2$.

\vspace{0.2cm}

Therefore, Corollary $5.3$ is a version of Colored Tverberg theorem of \v{Z}ivaljevi\'{c} and Vre\'{c}ica (Theorem $1.8$) that reduces the possible Tverberg partitions $\{\sigma_1, \ldots , \sigma_r\}$ such that $f(\sigma_1) \cap \cdots \cap f(\sigma_r) \neq \emptyset$.

\end{example}

\end{obs}

\vspace{0.2cm}

In particular, taking $m = d+1$, $l_1 = \cdots = l_d = r-1$ and $l_{d+1} =r$ we have the following consequence of Theorem $5.1$.

\vspace{0.2cm}

\begin{cor}

Let $d \ge 1$ be an integer and let $r=p^{n} \ge 2$ be a prime power. For every continuous map $f: \Delta \rightarrow \mathbb{R}^{d}$ and every coloring $(C_1, \ldots , C_{d+1})$ of the vertex set of the simplex $\Delta$ by $d+1$ colors, with each color of size at least $2r-1$, there exist $r$ pairwise disjoint faces $\sigma_1, \ldots , \sigma_r$ of $\Delta$  satisfying \begin{eqnarray} {\rm (i)} \hspace{0.2cm} | C_i \cap \sigma_j| \le 1 \mbox{,} \hspace{0.1cm} \mbox{ for every } i \in [d+1], \hspace{0.1cm} j \in [r], \hspace{0.1cm} \mbox{that is, } \sigma_1, \ldots, \sigma_r \hspace{2.1cm} \nonumber \\ \mbox{are rainbow faces} ;\nonumber \hspace{10cm} \\\hspace{0.1cm} {\rm (ii)} \hspace{0.2cm} f(\sigma_1) \cap \cdots \cap f(\sigma_r) \neq \emptyset \mbox{; and} \hspace{8cm}\nonumber \\ {\rm (iii)} \hspace{0.2cm} | (\sigma_1 \cup \cdots \cup \sigma_r) \cap C_i| \le r-1 \mbox{,} \hspace{0.1cm} \mbox{ for every } i \in [d]. \hspace{0.1cm}  \hspace{4.3cm}\nonumber \end{eqnarray}

\end{cor}

\vspace{0.2cm}

\textit{Proof.} Just apply Theorem $5.1$ to the case that $m=d+1$, $l_1 = \cdots = l_d = r-1$ and $l_{d+1} =r$.

\hspace{12.3cm} $\square$

\vspace{0.2cm}

Note that constraint {\rm (iii)} does not include the condition: \begin{eqnarray} | (\sigma_1 \cup \cdots \cup \sigma_r) \cap C_{d+1}| \le r-1 .  \hspace{0.3cm}   \end{eqnarray}

\vspace{0.2cm}

If we add the assumption $(1)$ to constraint {\rm (iii)} and if we use the same tools as the proof of Theorem $5.1$, the answer to this problem is inconclusive, so we have the following unsolved problem:

\vspace{0.2cm}

\begin{prob}

Let $d \ge 1$ be an integer and let $r =p^{n}\ge 2$ be a prime power. For every continuous map $f: \Delta \rightarrow \mathbb{R}^{d}$ and every coloring $(C_1, \ldots , C_{d+1})$ of the vertex set of the simplex $\Delta$ by $d+1$ colors, with each color of size at least $2r-1$, there exist $r$ pairwise disjoint faces $\sigma_1, \ldots , \sigma_r$ of $\Delta$  satisfying: \begin{eqnarray} {\rm (i)} \hspace{0.2cm} | C_i \cap \sigma_j| \le 1 \mbox{,}\hspace{0.1cm} \mbox{ for every } i \in [d+1], \hspace{0.1cm} j \in [r], \hspace{0.1cm} \mbox{that is, } \sigma_1, \ldots, \sigma_r \hspace{2.1cm} \nonumber \\ \mbox{are rainbow faces} ;\nonumber \hspace{10cm} \\  {\rm (ii)} \hspace{0.2cm} f(\sigma_1) \cap \cdots \cap f(\sigma_r) \neq \emptyset  \mbox{; and} \hspace{8cm}\nonumber \\ {\rm (iii)} \hspace{0.2cm} | (\sigma_1 \cup \cdots \cup \sigma_r) \cap C_i| \le r-1 \mbox{,} \hspace{0.1cm} \mbox{ for every } i \in [d+1]. \hspace{0.1cm} \hspace{3.7cm}\nonumber \end{eqnarray}

\end{prob}

\vspace{0.2cm}

However, if we add a color $C_{d+2}$ with a single vertex, the result follows, as in the follow result.

\vspace{0.2cm}

\begin{teo}

Let $d \ge 1$ be an integer and let $r=p^{n} \ge 2$ be a prime power. For every continuous map $f: \Delta \rightarrow \mathbb{R}^{d}$ and every coloring $(C_1, \ldots , C_{d+2})$ of the vertex set of the simplex $\Delta$ by $d+2$ colors, with colors $C_1, \ldots, C_{d+1}$ of size at least $2r-1$ and such that the color $C_{d+2}$ has a single vertex, there exist $r$ pairwise disjoint faces $\sigma_1, \ldots , \sigma_r$ of $\Delta$  satisfying: \begin{eqnarray} {\rm (i)} \hspace{0.2cm} | C_i \cap \sigma_j| \le 1 \mbox{,} \hspace{0.1cm} \mbox{ for every } i \in [d+2], \hspace{0.1cm} j \in [r], \hspace{0.1cm} \mbox{that is, } \sigma_1, \ldots, \sigma_r \hspace{2.1cm} \nonumber \\ \mbox{are rainbow faces} ;\nonumber \hspace{10cm} \\ {\rm (ii)} \hspace{0.2cm} f(\sigma_1) \cap \cdots \cap f(\sigma_r) \neq \emptyset \mbox{; and} \hspace{8cm} \nonumber \\ {\rm (iii)} \hspace{0.2cm} | (\sigma_1 \cup \cdots \cup \sigma_r) \cap C_i| \le r-1 \mbox{,} \hspace{0.1cm} \mbox{ for every } i \in [d+2]. \hspace{3.8cm} \nonumber \end{eqnarray}

\end{teo}

\vspace{0.2cm}

\textit{Proof.} Without loss of generality, we can assume that $|C_1| = \cdots = |C_{d+1}| = 2r-1$. Then, by Theorem $4.1$, it is sufficient to show that there is no $\mathfrak S_r$-equivariant map: \begin{eqnarray} \left(\Delta_{2r-1,r}^{\le r -2} \right)^{\ast (d+1)}  \ast \Delta_{1,r}\longrightarrow S \left(W_{r}^{\oplus (d+1)}\right). \nonumber \end{eqnarray}

\vspace{0.2cm}

Just as in Theorem $5.1$, we have a subgroup $(\mathbb{Z}_p)^{n} \cong \mbox{Im} (\varphi) \le \mathfrak S_r$. Hence, to prove the non-existence of $\mathfrak S_r$-equivariant map, it suffices to prove the non-existence of   $(\mathbb{Z}_p)^{n}$-equivariant map \begin{eqnarray} \left(\Delta_{2r-1,r}^{\le r -2} \right)^{\ast (d+1)}  \ast \Delta_{1,r}\longrightarrow S \left(W_{r}^{\oplus (d+1)}\right). \nonumber \end{eqnarray} 

\vspace{0.2cm}

Let us suppose that there is a such $(\mathbb{Z}_p)^{n}$-equivariant map. Thus, by Proposition $2.4(i)$, we have that: \begin{eqnarray} i\left( \left(\Delta_{2r-1,r}^{\le r -2}\right)^{\ast (d+1)}  \ast \Delta_{1,r}\right) \le i\left( S \left(W_{r}^{\oplus (d+1)}\right) \right). \nonumber \end{eqnarray}

\vspace{0.2cm}

It follows from Proposition $2.4 (iv)$ that  \begin{eqnarray} i\left( S \left(W_{r}^{\oplus (d+1)}\right) \right) = i\left( S^{(r-1) (d+1)-1}\right) = (r-1) (d+1).\nonumber \end{eqnarray}

\vspace{0.2cm}

On the other hand, note that $\Delta_{2r-1,r}^{\le r -2}$ is $(r-3)$-connected, for any $i=1, \ldots , d+1$.

\vspace{0.2cm}

Indeed, $\Delta_{2r-1,r}$ is $(r-2)$-connected (by Corollary $3.8$) so, it is $(r-3)$-connected. Therefore, as a result of Theorem $3.6$, we have that $\Delta_{2r-1,r}^{\le l_i -1}$ is $(r-3)$-connected.

\vspace{0.2cm}

The chessboard complex $\Delta_{1,r}$ is $(-1)$-connected, in consequence of Theorem $3.7$.

\vspace{0.2cm}

Due to Theorem $3.3$, we have that: \begin{eqnarray} \mbox{conn} \left ( \left(\Delta_{2r-1,r}^{\le r -2}\right)^{\ast (d+1)}  \ast \Delta_{1,r}\right) \ge  (d+1) (r-3) +(-1) + 2 (d+1) = \nonumber \\ (r-1) (d+1) -1. \hspace{9cm} \nonumber \end{eqnarray}

\vspace{0.2cm}

Therefore, accordingly Theorem $3.5$,  we have that \begin{eqnarray} i \left( \left(\Delta_{2r-1,r}^{\le r -2} \right)^{\ast (d+1)}  \ast \Delta_{1,r}\right) \ge  [(r-1) (d+1)-1] +2 =  (r-1) (d+1)+1. \nonumber \end{eqnarray}

\vspace{0.2cm}

Then, $i\left( \left(\Delta_{2r-1,r}^{\le r -2} \right)^{\ast (d+1)}  \ast \Delta_{1,r}\right) > (r-1) (d+1) =i \left(S\left(W_{r}^{\oplus(d+1)}\right)\right)$. But, it contradicts \begin{eqnarray} i\left( \left(\Delta_{2r-1,r}^{\le r -2} \right)^{\ast (d+1)} \ast \Delta_{1,r}\right) \le i \left(S\left(W_{r}^{\oplus(d+1)}\right)\right). \nonumber \end{eqnarray}

\vspace{0.2cm}

Thereby, is no $(\mathbb{Z}_p)^n$-equivariant map \begin{eqnarray} \left(\Delta_{2r-1,r}^{\le r -2} \right)^{\ast (d+1)} \ast \Delta_{1,r} \longrightarrow S \left(W_{r}^{\oplus (d+1)}\right). \nonumber \end{eqnarray}

\hspace{12.3cm} $\square$

\vspace{0.2cm}

The next results are variations of the Problem $1.9$ and their proofs are based on the same techniques presented previously.

\vspace{0.2cm}

\begin{teo}

Let $d \ge 1$ be an integer and let $r =p^{n}\ge 3$ be a prime power. For every continuous map $f: \Delta \rightarrow \mathbb{R}^{d}$ and every coloring $(C_1, \ldots , C_{d+1})$ of the vertex set of the simplex $\Delta$ by $d+1$ colors, with colors $C_1, \ldots , C_d$ of size at least $2r-4$ and color $C_{d+1}$ of size least $2r-1$, there exist $r$ pairwise disjoint faces $\sigma_1, \ldots , \sigma_r$ of $\Delta$  satisfying: \begin{eqnarray} {\rm (i)} \hspace{0.2cm} | C_i \cap \sigma_j| \le 1 \mbox{,} \hspace{0.1cm} \mbox{ for every } i \in [d+1], \hspace{0.1cm} j \in [r], \hspace{0.1cm} \mbox{that is, } \sigma_1, \ldots, \sigma_r \hspace{2.1cm} \nonumber \\ \mbox{are rainbow faces} ;\nonumber \hspace{10cm} \\  \hspace{0.1cm} {\rm (ii)} \hspace{0.2cm} f(\sigma_1) \cap \cdots \cap f(\sigma_r) \neq \emptyset  \mbox{; and} \hspace{7.9cm}\nonumber \\ {\rm (iii)} \hspace{0.2cm} | (\sigma_1 \cup \cdots \cup \sigma_r) \cap C_i| \le r-1 \mbox{,} \hspace{0.1cm} \mbox{ for every } i \in [d]. \hspace{4.4cm} \nonumber \end{eqnarray}

\end{teo}

\vspace{0.2cm}

\textit{Proof.} Without loss of generality, we can assume that $|C_1| = \cdots = |C_{d}| = 2r-4$ and $|C_{d+1}|= 2r-1$. Then, by Theorem $4.1$ it is sufficient to show that there is no $\mathfrak S_r$-equivariant map: \begin{eqnarray} \left(\Delta_{2r-4,r}^{\le r -2}\right)^{\ast d}  \ast \Delta_{2r-1,r}\longrightarrow S \left(W_{r}^{\oplus (d+1)}\right). \nonumber \end{eqnarray}

\vspace{0.2cm}

Just as in Theorem $5.1$, we have a subgroup $(\mathbb{Z}_p)^{n} \cong \mbox{Im} (\varphi) \le \mathfrak S_r$. Hence, to prove the non-existence of an $\mathfrak S_r$-equivariant map it suffices to prove the non-existence of a  $(\mathbb{Z}_p)^{n}$-equivariant map:  \begin{eqnarray} \left(\Delta_{2r-4,r}^{\le r -2} \right)^{\ast d}  \ast \Delta_{2r-1,r}\longrightarrow S \left(W_{r}^{\oplus (d+1)}\right).\nonumber \end{eqnarray}

\vspace{0.2cm}

Let us suppose that there is a such $(\mathbb{Z}_p)^{n}$-equivariant map. Thus, by Proposition $2.4 (i)$, we have that: \begin{eqnarray} i\left( \left(\Delta_{2r-4,r}^{\le r -2}\right)^{\ast d} \ast \Delta_{2r-1,r}\right) \le i\left( S \left(W_{r}^{\oplus (d+1)}\right) \right). \nonumber \end{eqnarray}

\vspace{0.2cm}

It follows from Proposition $2.4 (iv)$ that  \begin{eqnarray} i\left( S \left(W_{r}^{\oplus (d+1)}\right) \right) = i\left( S^{(r-1) (d+1)-1}\right) = (r-1) (d+1). \nonumber \end{eqnarray}

\vspace{0.2cm}

On the other hand, note that $\Delta_{2r-4,r}^{\le r -2}$ is $(r-3)$-connected, for any $i=1, \ldots , d$.

\vspace{0.2cm}

Indeed, $\Delta_{2r-4,r}$ is $(r-3)$-connected (by Theorem $3.7$) so, it is $(r-3)$-connected. Therefore, as a result of Theorem $3.6$, we have that $\Delta_{2r-4,r}^{\le r-2}$ is $(r-3)$-connected.

\vspace{0.2cm}

The chessboard complex $\Delta_{2r-1,r}$ is $(r-2)$-connected, in consequence of Theorem $3.7$.

\vspace{0.2cm}

Due to Theorem $3.3$, we have that: \begin{eqnarray} \mbox{conn} \left( \left(\Delta_{2r-4,r}^{\le r -2} \right)^{\ast d} \ast \Delta_{2r-1,r}\right) \ge  (r-3) d + (r-2) +2d = \nonumber \\ (r-1) (d+1) -1. \hspace{7cm}  \nonumber \end{eqnarray}

\vspace{0.2cm}

Therefore, accordingly Theorem $3.5$,  we have that: \begin{eqnarray} i \left( \left(\Delta_{2r-4,r}^{\le r -2}\right)^{\ast d} \ast \Delta_{2r-1,r}\right) \ge  [(r-1) (d+1)-1] +2 = (r-1) (d+1)+1. \nonumber \end{eqnarray}

\vspace{0.2cm}

Then, $i \left( \left(\Delta_{2r-4,r}^{\le r -2} \right)^{\ast d} \ast \Delta_{2r-1,r}\right) > (r-1) (d+1) =i \left(S\left(W_{r}^{\oplus(d+1)}\right)\right)$. But, it contradicts \begin{eqnarray} i\left( \left(\Delta_{2r-4,r}^{\le r -2} \right)^{\ast d}  \ast \Delta_{2r-1,r}\right) \le i \left(S\left(W_{r}^{\oplus(d+1)}\right)\right) \nonumber. \end{eqnarray} Thereby, there is no $(\mathbb{Z}_p)^n$-equivariant map \begin{eqnarray} \left(\Delta_{2r-4,r}^{\le r -2} \right)^{\ast d} \ast \Delta_{2r-1,r} \longrightarrow S \left(W_{r}^{\oplus (d+1)}\right). \nonumber \end{eqnarray}

\hspace{12.3cm} $\square$

\vspace{0.2cm}

\begin{teo}

Let $d \ge 1$ be an integer, and let $r =p^{n} \ge 3$ be a prime power. For every continuous map $f: \Delta \rightarrow \mathbb{R}^{d}$ and every coloring $(C_1, \ldots , C_{d+2})$ of the vertex set of the simplex $\Delta$ by $d+2$ colors, with colors $C_1, \ldots , C_{d+1}$ of size at least $2r-4$ and such that the color $C_{d+1}$ has a single vertex, there exist $r$ pairwise disjoint faces $\sigma_1, \ldots , \sigma_r$ of $\Delta$  satisfying: \begin{eqnarray} {\rm (i)} \hspace{0.2cm} | C_i \cap \sigma_j| \le 1 \mbox{,} \hspace{0.1cm} \mbox{ for every } i \in [d+2], \hspace{0.1cm} j \in [r], \hspace{0.1cm} \mbox{that is, } \sigma_1, \ldots, \sigma_r \hspace{5.1cm} \nonumber \\ \mbox{are rainbow faces} ;\nonumber \hspace{13cm} \\   {\rm (ii)} \hspace{0.2cm} f(\sigma_1) \cap \cdots \cap f(\sigma_r) \neq \emptyset\mbox{; and} \hspace{11cm}\nonumber\\ {\rm (iii)} | (\sigma_1 \cup \cdots \cup \sigma_r) \cap C_i| \le r-1 \mbox{,} \hspace{0.1cm} \mbox{ for every } i \in [d+2]. \hspace{0.1cm}   \hspace{6.8cm}\nonumber\end{eqnarray}

\end{teo}

\vspace{0.2cm}

\textit{Proof.} Without loss of generality, we can assume that $|C_1| = \cdots = |C_{d+1}| = 2r-4$. Then, by Theorem $4.1$ it is sufficient to show that there is no $\mathfrak S_r$-equivariant map: \begin{eqnarray} \left(\Delta_{2r-4,r}^{\le r -2} \right)^{\ast (d+1)} \ast \Delta_{1,r}\longrightarrow S \left(W_{r}^{\oplus (d+1)}\right). \nonumber \end{eqnarray}

\vspace{0.2cm}

Just as in Theorem $5.1$, we have a subgroup $(\mathbb{Z}_p)^{n} \cong \mbox{Im} (\varphi) \le \mathfrak S_r$. Hence, to prove the non-existence of a $\mathfrak S_r$-equivariant map, it suffices to prove the non-existence of a  $(\mathbb{Z}_p)^{n}$-equivariant map: \begin{eqnarray} \left( \Delta_{2r-4,r}^{\le r -2} \right)^{ \ast (d+1)} \ast \Delta_{1,r}\longrightarrow S \left(W_{r}^{\oplus (d+1)}\right). \nonumber \end{eqnarray}

\vspace{0.2cm}

Let us suppose that there is a such $(\mathbb{Z}_p)^{n}$-equivariant map. Thus, by Proposition $2.4 (i)$, we have that:  \begin{eqnarray} i\left( \left(\Delta_{2r-4,r}^{\le r -2}\right)^{ \ast (d+1)}  \ast \Delta_{1,r}\right) \le i\left( S \left(W_{r}^{\oplus (d+1)}\right) \right).  \nonumber \end{eqnarray}

\vspace{0.2cm}

It follows from Proposition $2.4 (iv)$ that  \begin{eqnarray} i\left( S \left(W_{r}^{\oplus (d+1)}\right) \right) = i\left( S^{(r-1) (d+1)-1}\right) = (r-1) (d+1). \nonumber\end{eqnarray}

\vspace{0.2cm}

On the other hand, note that $\Delta_{2r-4,r}^{\le r -2}$ is $(r-3)$-connected, for any $i=1, \cdots , d+1$.

\vspace{0.2cm}

Indeed, $\Delta_{2r-4,r}$ is $(r-3)$-connected (by Theorem $3.7$). Therefore, as a result of Theorem $3.6$, we have that $\Delta_{2r-4,r}^{\le r-2}$ is $(r-3)$-connected.

\vspace{0.2cm}

The chessboard complex $\Delta_{1,r}$ is $(-1)$-connected, in consequence of Theorem $3.7$.

\vspace{0.2cm}

Due to Theorem $3.3$, we have that: \begin{eqnarray} \mbox{conn} \left( \left(\Delta_{2r-4,r}^{\le r -2}\right)^{\ast (d+1)} \ast \Delta_{1,r}\right) \ge  (d+1) (r-3) +(-1) + 2 (d+1) = \nonumber \\ (r-1) (d+1) -1. \hspace{8.7cm} \nonumber \end{eqnarray}

\vspace{0.2cm}

Therefore, accordingly Theorem $3.5$,  we have that: \begin{eqnarray} i \left( \left(\Delta_{2r-1,r}^{\le r -2} \right)^{\ast (d+1)} \ast \Delta_{1,r}\right) \ge  [(r-1) (d+1)-1] +2 = (r-1) (d+1)+1. \nonumber \end{eqnarray}

\vspace{0.2cm}

Then, $i\left( \left(\Delta_{2r-4,r}^{\le r -2} \right)^{\ast (d+1)} \ast \Delta_{1,r}\right) > (r-1) (d+1) =i \left(S\left(W_{r}^{\oplus(d+1)}\right) \right)$. But, it contradicts \begin{eqnarray} i\left( \left(\Delta_{2r-4,r}^{\le r -2} \right)^{\ast (d+1)} \ast \Delta_{1,r}\right) \le i \left(S\left(W_{r}^{\oplus(d+1)}\right) \right) \nonumber. \end{eqnarray} Thereby, there is no $(\mathbb{Z}_p)^n$-equivariant map \begin{eqnarray} \left(\Delta_{2r-4,r}^{\le r -2}\right)^{\ast (d+1)} \ast \Delta_{1,r} \longrightarrow S \left(W_{r}^{\oplus (d+1)}\right). \nonumber \end{eqnarray}

\hspace{12.3cm} $\square$

\vspace{0.2cm}

\end{teo}

\end{document}